# GAUSSIAN FLUCTUATIONS FOR NON-HERMITIAN RANDOM MATRIX ENSEMBLES

By B. Rider and Jack W. Silverstein

*University of Colorado at Boulder and North Carolina State University*

Consider an ensemble of $N \times N$ non-Hermitian matrices in which all entries are independent identically distributed complex random variables of mean zero and absolute mean-square one. If the entry distributions also possess bounded densities and finite $(4 + \varepsilon)$ moments, then Z. D. Bai [*Ann. Probab.* **25** (1997) 494–529] has shown the ensemble to satisfy the circular law: after scaling by a factor of $1/\sqrt{N}$ and letting $N \to \infty$, the empirical measure of the eigenvalues converges weakly to the uniform measure on the unit disk in the complex plane. In this note, we investigate fluctuations from the circular law in a more restrictive class of non-Hermitian matrices for which higher moments of the entries obey a growth condition. The main result is a central limit theorem for linear statistics of type $X_N(f) = \sum_{k=1}^{N} f(\lambda_k)$ where $\lambda_1, \lambda_2, \ldots, \lambda_N$ denote the ensemble eigenvalues and the test function $f$ is analytic on an appropriate domain. The proof is inspired by Bai and Silverstein [*Ann. Probab.* **32** (2004) 533–605], where the analogous result for random sample covariance matrices is established.

**1. Introduction.** Eigenvalues of random Hermitian matrices in the limit of large dimension have been widely studied on account of their relevance to multivariate statistics and theoretical physics. The study of non-Hermitian random matrix ensembles remains undeveloped by comparison. While considered as early as the mid-1960s by Ginibre [15], only in the last several years has the non-Hermitian case been picked up with serious interest, notably in the physics community ([10, 11] and [24] offer guides to that literature).

Regardless of symmetry considerations, the most basic object of study in any random matrix ensemble is the density of states, or the limiting spectral distribution of the eigenvalues. In the setting considered here, this may be described as follows. Consider an $N \times N$ matrix in which all entries are independent identically distributed and standardized to have mean zero and mean-square one. If this matrix is normalized by a factor of $1/\sqrt{N}$ and the limit $N \to \infty$ taken, then (modulo, perhaps, additional technical assumptions) it is anticipated that the counting measure of eigenvalues will converge weakly to the uniform measure on the unit disk in the complex plane $\mathbb{C}$. This is referred to as the *circular law*.

In the fundamental work [15], the cases in which all entries are either independent identically distributed complex or quaternion Gaussians are shown to be integrable. That is, the full joint density of the $N$ eigenvalues is computed. Employing this explicit eigenvalue density, the first proper proof of the circular law for the complex Gaussian ensemble (Ginibre ensemble) may be found in [19]; see also [9] which contains a weaker form of the circular law for ensembles of real Gaussian entries. Using an important idea of Girko [14], Bai [1] has proved the most general result to date: the circular law holds for ensembles of independent complex entries such that the common entry distribution possesses a bounded density and finite $(4+\varepsilon)$ moment.

This paper is concerned with fluctuations about the circular law. We focus on linear spectral statistics of the form

$$(1.1) \qquad X_N(f) = \sum_{k=1}^{N} f(\lambda_k),$$

where $\lambda_1, \lambda_2, \ldots, \lambda_N$ denote the ensemble eigenvalues. Note that $X_N(f) = N \int f(\lambda) \mu_N(d\lambda)$ where $\mu_N(d\lambda)$ is the empirical spectral distribution. For ensembles of Hermitian matrices, there is an enormous body of work on Gaussian distributional limits of $X_N(f)$ as $N \to \infty$, for which a variety of methods has been developed. We cite [3, 7, 16, 20, 23, 27] and [28] as a representative list. We should also mention a different class of ensembles for which similar investigations have been carried out. These are the classical groups: $N \times N$ unitary, orthogonal or symplectic matrices distributed according to Haar measure. While certain techniques may be imported from the Hermitian setting, these ensembles have their own structure and their analysis requires a separate set of tools; see, for example, [8, 18, 21, 22] and [30]. One common theme in all of these results is how they manifest the rigidity of the point processes formed by random matrix eigenvalues. For smooth test functions $f$, the linear statistic need only be centered to produce a central limit theorem, no normalization is required. [Even for discontinuous test functions, the growth of the variance of $X_N(f)$ is much less than that for a set of $N$ independent particles.]

Thus far, fluctuation results for linear statistics of non-Hermitian random matrices have been relegated to Ginibre's complex ensemble, where the explicit joint density of eigenvalues may be exploited. In [10], advantage is taken of the particularly simple form of the density of the eigenvalue moduli in the Ginibre ensemble and a central limit theorem is described for statistics $X_N(f) - E[X_N(f)]$, where $f(\lambda) = f(|\lambda|)$ and $f$ is differentiable. [Again the variance is $O(1)$ as $N \to \infty$.] Reference [26] extends this picture to discontinuous statistics of the eigenvalue moduli and also proves a collection of results for the case where $f$ depends only on the eigenvalue phase. Here, we are able to go beyond the Gaussian (Ginibre) setup and establish a central limit theorem for spectral statistics (1.1) in a class of non-Hermitian random matrices built out of more general entry distributions. The proof is inspired by the work of the second author and Z. D. Bai [3]. The method is fairly robust in terms of the underlying distribution of the matrix entries but does require the analyticity of the test function $f$ in $X_N(f)$.

The particular ensemble we consider, denoted by $M$, may be described as $1/\sqrt{N} \times$ a base $N \times N$ matrix $\tilde{M}$. The latter is comprised of independent identically distributed complex entries $\tilde{m}_{ij}$ which are normalized in the usual way,

$$E[\tilde{m}_{11}] = 0 \quad \text{and} \quad E[|\tilde{m}_{11}|^2] = 1$$

and which also satisfy

(1.2)
(i) $E[\tilde{m}_{11}^2] = 0$,
(ii) $E[|\tilde{m}_{11}|^k] \leq k^{\alpha k}$ for $k > 2$ and some $\alpha > 0$,
(iii) $\Re(\tilde{m}_{11})$ and $\Im(\tilde{m}_{11})$ possess a bounded joint density.

Condition (1.2.i) is a Gaussian-like assumption; it is analogous to requirements (ii) and (iii) in Theorem 1.1 of [3]. The moment condition (1.2.ii) supplies a priori control on the magnitude of the eigenvalues. In particular, (1.2.ii) implies that with probability one as $N \to \infty$, the spectral radius of $M$ converges to 1 and the spectral norm of $M$ converges to 2, see [12] and [13]. While it is more natural to consider the spectral radius in the non-Hermitian setting, it will be the spectral norm which figures more prominently below. Finally, the smoothness condition (1.2.iii) is of a technical nature. The same assumption was put to similar use in the proof of the circular law in [1].

With the ensemble defined we may state our first result:

THEOREM 1.1. *Let the independent entries of the non-Hermitian matrix $M$ satisfy the conditions outlined in* (1.2). *Consider test functions $f_1, f_2, \ldots, f_k$, analytic in a neighborhood of the disk $|z| \leq 4$ and otherwise bounded. Then as $N \to \infty$, the vector*

(1.3) $\quad (X_N(f_1) - Nf_1(0), X_N(f_2) - Nf_2(0), \ldots, X_N(f_k) - Nf_k(0))$

*converges in distribution to the mean-zero multivariate Gaussian* $(X(f_1), X(f_2), \ldots, X(f_k))$ *with covariances*

$$(1.4) \qquad E[X(f_\ell)\overline{X(f_m)}] = \frac{1}{\pi} \int_{\mathbb{U}} \frac{d}{dz} f_\ell(z) \overline{\frac{d}{dz} f_m(z)} \, d^2 z,$$

*in which* $\mathbb{U}$ *is the unit disk and* $d^2 z = d\Re(z) \, d\Im(z)$.

Note that $X_N(f)$ has been centered about the asymptotic mean. For large $N$, $E[X_N(f)]$ is well approximated by $N \times$ the average of $f$ against the circular law and

$$\frac{1}{\pi} \int_{\mathbb{U}} f(z) \, d^2 z = f(0)$$

if $f$ is analytic. The limit distribution of $X_N(f) - E[X_N(f)]$ is the same, as will become apparent in the course of the proof. The assumed analyticity of the test functions is a by-product of the strategy of proof: our main accomplishment in this note is actually a central limit theorem for the resolvent of the random matrix $M$. Central limit theorems for analytic statistics are then inferred by way of the Cauchy integral formula. More to the point, for any domain $\mathbb{D} \subset \mathbb{C}$ on which $f$ is analytic,

$$(1.5) \qquad X_N(f) - Nf(0) = \frac{1}{2\pi i} \int_{\partial \mathbb{D}} f(z) \{\operatorname{tr}(z - M)^{-1} - Nz^{-1}\} \, dz$$

holds on the event that all of $\lambda_1, \ldots, \lambda_N$ are contained in $\mathbb{D}$. Assuming that the probability of the complement of that event decays sufficiently fast [which is the content of estimate (2.1) below], Theorem 1.1 is a corollary of the next result.

THEOREM 1.2. *View the centered resolvent*

$$\mathfrak{G}_N(z) = \operatorname{tr}(z - M)^{-1} - Nz^{-1}$$

*as a process in* $z \in \mathcal{C}$, *any fixed contour lying exterior to the disk of radius four. Then the family* $\{\mathfrak{G}_N(z)\}$ *is tight in the space of continuous functions on* $\mathcal{C}$ *and converges weakly to the complex Gaussian process* $\mathfrak{G}(z)$ *defined by*

$$(1.6) \qquad E[\mathfrak{G}(z)\overline{\mathfrak{G}(w)}] = (1 - z\bar{w})^{-2}$$

*and* $E[\mathfrak{G}(z)] = 0$.

To pass the central limit theorem for $\mathfrak{G}_N(z)$ onto $X_N(f)$ requires only that we note the continuity of the map $\mathfrak{G}_N \to \int_{\mathcal{C}} f(z) \mathfrak{G}_N(z) \, dz$ from the space of continuous functions on $\mathcal{C}$ into $\mathbb{C}$. As for connecting the covariance structures (1.4) and (1.6), choose functions $f$ and $g$ analytic in a neighborhood of $\mathbb{U}$. Also, choose a contour $\mathcal{C}$ lying within the region of analyticity of both functions, but enclosing $\mathbb{U}$. Then (1.4), applied to $X(f)$ and $X(g)$, may be rewritten as

$$
\begin{aligned}
(1.7) \quad & \frac{1}{\pi} \int_{\mathbb{U}} f'(\eta)\overline{g'(\eta)}\, d^2\eta \\
& = -\frac{1}{4\pi^2} \int_{\mathcal{C}}\int_{\mathcal{C}} f(\eta)\overline{g(\eta)} \left\{ \frac{1}{\pi} \int_{\mathbb{U}} \frac{d^2\eta}{(\eta-z)^2(\bar{\eta}-\bar{w})^2} \right\} dz\, d\bar{w}
\end{aligned}
$$

and the desired equality,

$$
(1.8) \qquad \frac{1}{\pi} \int_{\mathbb{U}} \frac{d^2\eta}{(\eta-z)^2(\bar{\eta}-\bar{w})^2} = (1-z\bar{w})^{-2},
$$

may be verified by expanding the integrand on the left-hand side.

We point out that the covariance of the limiting resolvent, $(1-z\bar{w})^{-2}$, equals $\pi \times$ the *Bergman kernel* for $\mathbb{U}$ under that map $z \to 1/z$ (for general background on the Bergman kernel, [4] is recommended). The Bergman kernel also appears in the covariance of the zero process of a conformally invariant Gaussian analytic function on the unit disk, recently studied in [25]. In fact, the result of Theorem 1.2 shows that the limiting resolvent $\mathfrak{G}(z)$ coincides, at least in the sense of finite-dimensional distributions, with the Gaussian analytic function

$$
(1.9) \qquad G(z) \equiv \sum_{k=1}^{\infty} \sqrt{k} Z_k \left(\frac{1}{z}\right)^{k+1},
$$

in which $\{Z_k\}$ are independent unit complex Gaussians and the prescription is valid for $|z| > 1$. In this vein, it is also interesting to compare Theorem 1.2 with the work of [8] on the logarithmic derivative of the characteristic polynomial of a random unitary matrix.

Finally, the fact that Theorem 1.2 is stated for $|z| > 4$ (and thus the analyticity in Theorem 1.1 required for $|z| \leq 4$), whereas the limiting covariance (1.6) and the appraisals (1.7) through (1.9) are on firm ground as soon as $z$ lies in the exterior of $\mathbb{U}$, demands explanation. Undoubtedly, $|z| > 1$ is the correct condition. In fact, for the complex Ginibre ensemble, the authors can prove that the corresponding resolvent satisfies Theorem 1.2 for all $|z| > 1$; this work will appear elsewhere. Further, central limit theorems of the type outlined in Theorem 1.1 are anticipated to hold for test functions which are just once continuously differentiable. The stronger assumptions made here for more general entry distributions are an artifact of the method. Experts on random matrix theory will appreciate the challenges in extending techniques from the study of Hermitian random matrices to the non-Hermitian setting.

Since the eigenvalues of $M$ can lie throughout the complex plane, the first step in any proof of Theorem 1.2 is to establish an a priori control of the resolvent. With our assumptions granting that the spectral radius $\equiv |\lambda_{\max}(M)|$ tends to one almost surely, it is natural to attempt to work on this event. That is, taking $|z|$ larger than $1+\varepsilon$ for some $\varepsilon > 0$, $\operatorname{tr}(z-M)^{-1}$ is well-behaved on the event that $|\lambda_{\max}(M)|$ is less than, say, $1+\varepsilon/2$. The idea would be to prove a central limit theorem for this truncated object, then to show that the resolvent times the indicator of $\{|\lambda_{\max}(M)| > 1+\varepsilon/2\}$ goes to zero in probability in the appropriate norm. However, the mechanics of our proof do not allow us to take full advantage of the known estimates of the spectral radius.

The basic idea behind the proof of Theorem 1.2 is to find a martingale structure in $\operatorname{tr}(z-M)^{-1}$ by successively conditioning on the sigma-fields generated by the columns of $M$. Given the desired structure, checking that the conditions of a standard martingale central limit theorem are satisfied requires estimates which compare the resolvent matrix of $M$ and that of various rank-one perturbations of $M$. An example of such a perturbation is the matrix formed by setting the $k$th column of $M$ to zero, denoted by $M_k$. The problem which arises is as follows. Having control of $|\lambda_{\max}(M)|$ does not imply that we have similar control over $|\lambda_{\max}(M_k)|$. On the other hand, the spectral norm of $M$ trivially bounds the spectral norm of $M_k$ for all $k$. Since assumption (1.2.ii) implies that the former lies under two with overwhelming probability, we can attain uniform control of $\operatorname{tr}(z-M)^{-1}$, $\operatorname{tr}(z-M_k)^{-1}$, etc., if we are willing to take $|z|$ larger than four.

The next section picks up this discussion, describes an overview of the proof of Theorem 1.2 and records various estimates needed therein. The proof itself then occupies Sections 3 and 4, along with a brief Appendix.

**2. Outline and preliminary facts.** Henceforth, we denote the resolvent matrix by

$$\Xi_N(z)^{-1} \equiv (z-M)^{-1}.$$

As indicated in the Introduction, a characteristic feature of the proof is the procedure by which the norm of the eigenvalues of $M$ is controlled. For this, we define

$$\Omega_N \equiv \{M : \|M\| < \kappa\},$$

in which $\|\cdot\|$ is the spectral norm and $\kappa$ is a number larger than two which we now fix. Note that for $|z| > 4$, there exists a constant $\mathcal{K} = \mathcal{K}(\kappa, z)$ such that

$$\|\Xi_N(z)^{-1}\| \leq \mathcal{K} \qquad \text{on the event } \Omega_N.$$

Also, from this point, whenever we write $|z| > 4$, it should be understood that $|z|$ exceeds four by some fixed $\varepsilon > 0$. Thus, the $\mathcal{K}$ in the above bound is also fixed from this point on. We reiterate that $\|M\| < \kappa$ implies that $\|M_k\| < \kappa$ for any $k$, with the same implication holding for further perturbations of $M$ (in which more than one column is removed, or set to zero). Such comparisons constitute the chief reason we are employing the spectral norm, rather than the more natural spectral radius, to control the support of the empirical eigenvalue distribution.

Since we essentially wish to work on $\Omega_N$, it is of course important that these events exhaust the space as $N \to \infty$. In this, we are fortunate that the following estimate has been obtained by Geman [13]: under the moment assumption (1.2.ii),

$$(2.1) \qquad P(\Omega_N^C) = o(N^{-\alpha})$$

for any $\alpha > 0$. It should be noted that while the proof in [13] is made for real-valued random entries, the extension to the complex case is immediate. Also, if the entries of $M$ are drawn from a single doubly-infinite array, [31] proves that $\|M\| \to 2$ with probability one as soon as $E|m_{11}|^4 < \infty$. However, this result does not come equipped with an estimate as sharp as (2.1).

The stage is now set to describe the proof of the central limit theorem for $\mathfrak{G}_N(z) = \operatorname{tr} \Xi_N(z)^{-1} - Nz^{-1}$. First, we have the decomposition

$$\begin{aligned}\mathfrak{G}_N(z) &= (\operatorname{tr} \Xi_N(z)^{-1} \mathbf{1}_{\Omega_N} - E[\operatorname{tr} \Xi_N(z)^{-1}, \Omega_N]) \\ &\quad + (\operatorname{tr} \Xi_N(z)^{-1} \mathbf{1}_{\Omega_N^C} - E[\operatorname{tr} \Xi_N(z)^{-1}, \Omega_N^C]) \\ &\quad + (E[\operatorname{tr} \Xi_N(z)^{-1}] - Nz^{-1}) \\ &\equiv \mathfrak{G}_N^0(z) + \mathfrak{G}_N^1(z) + \mathfrak{G}_N^2(z).\end{aligned}$$

The Gaussian limit comes out of the first term on the right-hand side, $\mathfrak{G}_N^0(z)$, which is the original resolvent, truncated and recentered. Let $E_k$ denote the conditional expectation with respect to the sigma-field generated by the first $k$ columns, $m_1, m_2, \ldots, m_k$ of $M$ and write

$$(2.2) \qquad \mathfrak{G}_N^0(z) = \sum_{k=1}^N (E_k - E_{k-1})[\operatorname{tr} \Xi_N(z)^{-1}, \Omega_N]$$

($E_0[\cdot] = E[\cdot]$). Clearly, $\{(E_k - E_{k-1})[\operatorname{tr} \Xi(z)^{-1}, \Omega_N]\}$ form a (bounded) martingale difference sequence and we recall the following standard result:

LEMMA 2.1 (Theorem 35.12 of [6]). *For each $N$, let $X_{N,1}, X_{N,2}, \ldots, X_{N,k_N}$ be a real-valued martingale difference sequence with respect to the filtration $\mathcal{F}_{N,k}$, each term having bounded second moments. Suppose that for any $\varepsilon > 0$ and a constant $\sigma^2$,*

(i) $\displaystyle \lim_{N \to \infty} \sum_{k=1}^{k_N} E[X_{N,k}^2, X_{N,k}| > \varepsilon] = 0,$

(ii) $\lim_{N \to \infty} P\left(\left|\sum_{k=1}^{k_N} E[X_{N,k}^2 | \mathcal{F}_{N,k}] - \sigma^2\right| > \varepsilon\right) = 0.$

Then as $N \to \infty$, the distribution of $\sum_{k=1}^{k_N} X_{N,k}$ converges weakly to a Gaussian with mean zero and variance $\sigma^2$.

Section 3 proves that appropriate linear combinations of $\mathfrak{G}_N^0(z)$ satisfy the conditions of Lemma 2.1. Along with a tightness estimate, this will show that over a fixed contour $\mathcal{C}$ exterior to the disk of radius four, $\mathfrak{G}_N^0(z)$ converges weakly to the mean zero Gaussian process with covariance $(1 - z\bar{w})^{-2}$ for $z, w \in \mathcal{C}$.

The proof of Theorem 1.2 is then reduced to showing that

(2.3) $$\sup_{z \in \mathcal{C}} |\mathfrak{G}_N^1(z)| \to 0$$

in probability and that

(2.4) $$\lim_{N \to \infty} \sup_{z \in \mathcal{C}} |\mathfrak{G}_N^2(z)| = 0$$

for the same given contour $\mathcal{C}$. The proof of (2.4) occupies Section 4; it is convenient to dispense with (2.3) immediately.

Since we have assumed that the densities of the individual entries of $M$ are continuous, it follows that $\sup_{z \in \mathcal{C}} |\operatorname{tr} \Xi_N(z)^{-1}|$ is finite for all $N$ with probability one. (With probability one, the spectrum of $M$ does not intersect any predetermined contour.) Therefore,

$$P\left(\sup_{z \in \mathcal{C}} |\operatorname{tr} \Xi_N(z)^{-1}| \mathbf{1}_{\Omega_N^C} \geq \varepsilon\right) \leq P(\Omega_N^C),$$

which goes to zero as $N \to \infty$, by (2.1). To complete the verification of (2.3), we must have

$$\lim_{N \to \infty} \sup_{z \in \mathcal{C}} E[|\operatorname{tr} \Xi_N(z)^{-1}|, \Omega_N^C] = 0.$$

This, however, is a consequence of Hölder's inequality, (2.1), and the following (rough) estimate, the proof of which is found in the Appendix:

LEMMA 2.2. *For any $p$, $1 \leq p < 2$ and any $z \in C$,*

$$E |\operatorname{tr} \Xi_N(z)^{-1}|^p \leq C(p) N^{3p/2+2}.$$

*The constant $C(p) \uparrow \infty$ as $p \uparrow 2$, but is independent of $z$. We note for later that the same bound holds if $M$ in $\Xi_N(z)^{-1} = (z - M)^{-1}$ is replaced by $M_k$.*

Finally, we conclude this section by recording the following estimate, of which we will make repeated use:

LEMMA 2.3 ([2], Lemma 2.7). *For a vector $X = (x_1, x_2, \ldots, x_N)^T$ of independent identically distributed (complex) entries $x_k$ with $E[x_1] = 0$ and $E[|x_1|^2] = 1$ and any $N \times N$ (complex) matrix $A$, it holds that*

$$E|X^*AX - \operatorname{tr} A|^p \leq C_1((E|x_1|^4 \operatorname{tr} A^*A)^{p/2} + E|x_1|^{2p} \operatorname{tr}(A^*A)^{p/2})$$

*for all $p \geq 2$ and a constant $C_1$ depending only on $p$.*

Even more useful will be the related bound: with $C_2$ being a new constant still depending only on $p$,

$$(2.5) \qquad E|X^*AX|^p \leq C_2 E|x_1|^{2p}((\operatorname{tr} A^*A)^{p/2} + |\operatorname{tr} A|^p).$$

This follows from the above lemma since $\operatorname{tr}(A^*A)^{p/2} \leq (\operatorname{tr} A^*A)^{p/2}$ for any matrix $A$.

**3. Convergence of the truncated process.** We establish the central limit theorem for the truncated (and recentered) resolvent $\mathfrak{G}_N^0$. Along with $M_k$ (or $M_{k\ell}$) denoting the matrix formed by setting the $k$th column (or both the $k$th and $\ell$th columns) of $M$ to zero, we introduce the additional shorthand

$$\Xi_{N,k}(z)^{-1} = (z - M_k)^{-1}, \qquad \Xi_{N,k\ell}(z)^{-1} = (z - M_{k\ell})^{-1}$$

(and so forth), for the corresponding resolvent matrices. Likewise, $\Omega_{N,k}$ and $\Omega_{N,k\ell}$ (and so forth) will denote the events on which $M_k$ and $M_{k\ell}$ have spectral norms bounded by $\kappa$. Note that $\mathcal{K}$ still serves as a bound for the spectral norm of $\Xi_{N,k}(z)^{-1}$ (resp. $\Xi_{N,k\ell}(z)^{-1}$) on the event $\Omega_{N,k}$ (resp. $\Omega_{N,k\ell}$).

3.1. *Finite-dimensional distributions.* By the Cramér–Wold device, identifying the limiting finite-dimensional distributions of

$$\mathfrak{G}_N^0(z) = \sum_{k=1}^N (E_k - E_{k-1})[\operatorname{tr} \Xi_N(z)^{-1}, \Omega_N] \equiv \sum_{k=1}^N Z_{N,k}(z)$$

as multivariate Gaussian amounts to verifying the conditions of Lemma 2.1 for martingale-difference sums of the form

$$\mathcal{M}_N = \sum_{k=1}^N \left\{ \sum_{\ell=1}^K \alpha_\ell Z_{N,k}(z_\ell) + \beta_\ell \overline{Z_{N,k}(z_\ell)} \right\} \equiv \sum_{k=1}^N \mathcal{M}_{N,k}$$

for any fixed $K$, any choice of $z_\ell$'s on $\mathcal{C}$ and complex $(\alpha_\ell, \beta_\ell)$ such that the resulting object is real. As a first step, $\mathcal{M}_N$ is simplified through a series of reductions in each of its terms.

LEMMA 3.1. *As $N \to \infty$, $\mathcal{M}_N$ either converges in distribution or fails to along with a related martingale-difference sum $\widetilde{\mathcal{M}}_N$. The latter is defined by replacing each appearance of $Z_{N,k}(z)$ in $\mathcal{M}_{N,k}$ with*

$$\widetilde{Z}_{N,k}(z) = E_k[e_k^T \Xi_{N,k}(z)^{-2} m_k, \Omega_{N,k}] = e_k^T E_k[\Xi_{N,k}(z)^{-2}, \Omega_{N,k}] m_k,$$

*where $e_k$ is the canonical unit vector in $\mathbb{R}^N$. In other words, the difference of $\widetilde{\mathcal{M}}_N$ and $\mathcal{M}_N$ goes to zero in probability.*

The point of Lemma 3.1 is the introduction of independence in each term, particularly that of the vector $m_k$ and the matrix $E_k[\Xi_{N,k}(z)^{-2}, \Omega_{N,k}]$. With that established, the desired conclusion comes down to the following:

LEMMA 3.2. *The martingale-difference sequence $\widetilde{\mathcal{M}}_N$ satisfies the conditions of Lemma 2.1. In particular, the limiting covariance (1.6) in Theorem 1.2 is seen from the fact that*

$$(3.1) \qquad \sum_{k=1}^{N} E_{k-1}[\widetilde{\mathcal{M}}_{N,k}^2] \to \sum_{1 \leq i,j \leq K} \alpha_i \beta_j (1 - z_i \overline{z_j})^{-2}$$

*in probability as $N \to \infty$.*

PROOF OF LEMMA 3.1. Since $(E_k - E_{k-1})[\operatorname{tr} \Xi_{N,k}(z)^{-1}] = 0$,

$$(3.2) \qquad \begin{aligned} \sum_{k=1}^{N} Z_{N,k}(z) &= \sum_{k=1}^{N} (E_k - E_{k-1})[\operatorname{tr}(\Xi_N(z)^{-1} - \Xi_{N,k}(z)^{-1}), \Omega_N] \\ &\quad - \sum_{k=1}^{N} (E_k - E_{k-1})[\operatorname{tr} \Xi_{N,k}(z)^{-1}, \Omega_N^C]. \end{aligned}$$

By Hölder's inequality, estimate (2.1) of $P(\Omega_N^C)$ and Lemma 2.2, we may conclude that

$$E[|\operatorname{tr} \Xi_{N,k}(z)^{-1}|, \Omega_N^C] \leq C_0 N^{-2},$$

in which the constant $C_0$ is independent of the index $k$. That is, the second term on the right-hand side of (3.2) tends to zero in probability as $N \to \infty$ and thus may be neglected in questions of distributional convergence. For the first term, we bring in the usual resolvent identity along with the basic fact that

$$(3.3) \qquad (A + v e_k^T)^{-1} v = \frac{A^{-1} v}{(1 + e_k^T A^{-1} v)}$$

for any invertible $A$ and vector $v$. Applied in succession, this yields

$$\begin{aligned}
\operatorname{tr}(\Xi_N(z)^{-1} - \Xi_{N,k}(z)^{-1}) &= \operatorname{tr}(\Xi_N(z)^{-1}(m_k e_k^T)\Xi_{N,k}(z)^{-1}) \\
&= e_k^T(\Xi_{N,k}(z)^{-1}\Xi_N^{-1}(z))m_k \\
&= \frac{e_k^T \Xi_{N,k}(z)^{-2} m_k}{1 + e_k^T \Xi_{N,k}(z)^{-1} m_k} \\
&\equiv \delta_k(z) e_k^T \Xi_{N,k}(z)^{-2} m_k.
\end{aligned}$$
(3.4)

Next, presuming control on the norm of $\Xi_{N,k}(z)^{-1}$ and the fact that $m_k$ is mean-zero and normalized as in $E\|m_k\|_2^2 = \sum_{\ell=1}^N E|m_{\ell k}|^2 = 1$, it is anticipated that $e_k^T \Xi_{N,k}(z)^{-1} m_k$ tends to zero as $N \to \infty$. That is, we wish to replace $\delta_k(z)$ with one in each appearance of the sum (3.2), and conclude that

(3.5) $$\sum_{k=1}^N (E_k - E_{k-1})[e_k^T \Xi_{N,k}(z)^{-2} m_k - \delta_k(z) e_k^T \Xi_{N,k}(z)^{-2} m_k, \Omega_N] \to 0$$

in probability. Expanding a little, we find

$$\begin{aligned}
1 - \delta_k(z) = &-e_k^T \Xi_{N,k}(z)^{-1} m_k + (e_k^T \Xi_{N,k}(z)^{-1} m_k)^2 \\
&- (e_k^T \Xi_{N,k}(z)^{-1} m_k)^3 + \delta_k(z)(e_k^T \Xi_{N,k}(z)^{-1} m_k)^4.
\end{aligned}$$

The sum in question has been broken into four, the first three terms of which may be bounded in mean-square as follows. With $p = 1, 2$ or $3$,

$$E\left|\sum_{k=1}^N (E_k - E_{k-1})[(e_k^T \Xi_{N,k}(z)^{-1} m_k)^p (e_k^T \Xi_{N,k}(z)^{-2} m_k), \Omega_N]\right|^2$$

(3.6) $$\leq 4 \sum_{k=1}^N E[|e_k^T \Xi_{N,k}(z)^{-1} m_k|^{2p} |e_k^T \Xi_{N,k}(z)^{-2} m_k|^2, \Omega_{N,k}]$$

$$\leq 4 \sum_{k=1}^N \{E|m_k^* \mathbf{A}_{N,k}(z) m_k|^{2p} E|m_k^* \mathbf{B}_{N,k}(z) m_k|^2\}^{1/2},$$

in which $\mathbf{A}_{N,k}(z)$ and $\mathbf{B}_{N,k}(z)$ are the Hermitian matrices

$$\mathbf{A}_{N,k}(z) = (\Xi_{N,k}(z)^{-1})^*(e_k e_k^T)\Xi_{N,k}(z)^{-1}\mathbf{1}_{\Omega_{N,k}}$$

and

$$\mathbf{B}_{N,k}(z) = (\Xi_{N,k}(z)^{-2})^*(e_k e_k^T)\Xi_{N,k}(z)^{-2}\mathbf{1}_{\Omega_{N,k}}.$$

Having replaced $\Omega_N$ with $\Omega_{N,k} \supset \Omega_N$ in line two of (3.6), note that $m_k$ is independent of both $\mathbf{A}_{N,k}(z)$ and $\mathbf{B}_{N,k}(z)$. This makes way for an application of Lemma 2.3 [actually, its consequence, (2.5)]. We have

$$E|m_k^* \mathbf{A}_{N,k}(z) m_k|^{2p}$$
$$\leq C_1 N^{-2p} E[(\operatorname{tr} \mathbf{A}_{N,k}(z))^p]$$
(3.7)
$$= C_1 N^{-2p} E[(e_k^T (\Xi_{N,k}(z)^{-1} (\Xi_{N,k}(z)^{-1})^*) e_k)^p, \Omega_{N,k}]$$
$$\leq C_1 N^{-2p} E[\|\Xi_{N,k}(z)^{-1}\|^{2p}, \Omega_{N,k}]$$
$$\leq C_1 \mathcal{K}^{2p} N^{-2p},$$

where the finiteness of the constant $C_1$ relies on the assumption that the entries of $M$ have bounded moments. Since the expectation involving $\mathbf{B}_{N,k}(z)$ in (3.6) may clearly be estimated in the same way, it follows that the last line of that display is bounded above by a constant multiple of $N^{-p}$.

As for the last sum resulting from the expansion of (3.5), notice that on the set $\Omega_N$, the simple bound

$$|\delta_k(z) e_k^T \Xi_{N,k}(z)^{-2} m_k| \leq |\operatorname{tr} \Xi_N(z)^{-1}| + |\operatorname{tr} \Xi_{N,k}(z)^{-1}| \leq 2N\mathcal{K}$$

may be read off directly from the identity (3.4). It follows that

$$E\left|\sum_{k=1}^{N}(E_k - E_{k-1})[(e_k^T \Xi_{N,k}(z)^{-1} m_k)^4 \delta_k(z) e_k^T \Xi_{N,k}(z)^{-2} m_k, \Omega_N]\right|^2$$
(3.8)
$$\leq 16 N^2 \mathcal{K}^2 \sum_{k=1}^{N} E|m_k^* \mathbf{B}_{N,k}(z) m_k|^4$$
$$\leq C_2 N^{-1}.$$

This completes the verification of (3.5).

To produce the aforementioned form of $\widetilde{M}_{N,k}$ requires two more steps. First, we have

$$E\left|\sum_{k=1}^{N}(E_k - E_{k-1})[e_k^T \Xi_{N,k}^{-2}(z) m_k, \Omega_{N,k}]\right.$$
$$\left. - \sum_{k=1}^{N}(E_k - E_{k-1})[e_k^T \Xi_{N,k}^{-2}(z) m_k, \Omega_N]\right|$$
$$= E\left|\sum_{k=1}^{N}(E_k - E_{k-1})[e_k^T \Xi_{N,k}^{-2}(z) m_k, \Omega_{N,k} \cap \Omega_N^C]\right|$$
$$\leq 2\sum_{k=1}^{N}(E[|e_k^T \Xi_{N,k}(z)^{-2} m_k|^4, \Omega_{N,k}])^{1/4} (P(\Omega_N^C))^{3/4}$$

$$= 2 \sum_{k=1}^{N} (E|m_k^* \mathbf{B}_{N,k}(z) m_k|^2)^{1/4} (P(\Omega_N^C))^{3/4}$$

$$\leq C_3 \sqrt{N} (P(\Omega_N^C))^{3/4},$$

which goes to zero faster than any polynomial power of $N$ [recall (2.1)]. It remains to observe that

$$E_{k-1}[e_k^T \Xi_{N,k}(z)^{-1} m_k \Omega_{N,k}] = e_k^T E_{k-1}[\Xi_{N,k}(z)^{-1}, \Omega_{N,k}] E_{k-1}[m_k] = 0,$$

the independence being the point of replacing $\Omega_N$ with $\Omega_{N,k}$ termwise. The proof is thus complete. □

PROOF OF LEMMA 3.2. Verifying the first condition of the martingale central limit theorem (Lemma 2.1) is trivial. In the previous proof, we have already made repeated use of the fact that $E|\widetilde{Z}_{N,k}(z)|^4$ is less than $2E|m_k^* \mathbf{B}_{N,k}(z) m_k|^2 = O(N^{-2})$ for arbitrary $z \in \mathcal{C}$ and all $k$. Hence,

$$\sum_{k=1}^{N} E[\widetilde{\mathcal{M}}_{N,k}^2, |\widetilde{\mathcal{M}}_{N,k}| > \varepsilon] \leq \varepsilon^{-2} \sum_{k=1}^{N} E[\widetilde{\mathcal{M}}_{N,k}^4]$$

$$\leq \varepsilon^{-2} C_1 \sum_{k=1}^{N} \sum_{\ell=1}^{K} (|\alpha_\ell| + |\beta_\ell|)^4 E[|\widetilde{Z}_{N,k}(z_\ell)|^4]$$

tends to zero as $N \to \infty$.

Next, by the linearity of $\widetilde{\mathcal{M}}_{N,k}$ in the $\widetilde{Z}_{N,k}$'s, the covariance formula (3.1) follows from checking that

$$(3.9) \qquad \sum_{k=1}^{N} E_{k-1}[\widetilde{Z}_{N,k}(z) \widetilde{Z}_{N,k}(w)] = 0$$

and

$$(3.10) \qquad \sum_{k=1}^{N} E_{k-1}[\widetilde{Z}_{N,k}(z) \overline{\widetilde{Z}_{N,k}(w)}] \to \frac{1}{(1-z\bar{w})^2}$$

in probability as $N \to \infty$. For the equality, using the independence of $m_k$ and the $k$th row of

$$\mathcal{B}_{N,k}(z) = E_k[\Xi_{N,k}(z)^{-2}, \Omega_{N,k}],$$

we find that

$$\sum_{k=1}^{N} E_{k-1}[\widetilde{Z}_{N,k}(z) \widetilde{Z}_{N,k}(w)]$$

$$= \sum_{k=1}^{N} E_{k-1}[(e_k^T \mathcal{B}_{N,k}(z) m_k)(e_k^T \mathcal{B}_{N,k}(w) m_k)]$$

$$= \sum_{k=1}^{N} \sum_{1 \leq i,j \leq N} [\mathcal{B}_{N,k}(z)]_{ki} [\mathcal{B}_{N,k}(w)]_{kj} E[m_{ik} m_{jk}] = 0,$$

as both $E[m_{ij}]$ and $E[m_{ij}^2]$ are assumed to vanish. On the other hand, with conjugates present in the second half of (3.10), the computation becomes

$$\sum_{k=1}^{N} E_{k-1}[\widetilde{Z}_{N,k}(z) \overline{\widetilde{Z}_{N,k}(w)}]$$

(3.11)
$$= \frac{1}{N} \sum_{k=1}^{N} \mathrm{tr}(\mathcal{B}_{N,k}(w)^* (e_k e_k^T) \mathcal{B}_{N,k}(z))$$

$$= \frac{1}{N} \sum_{k=1}^{N} e_k^T (E_k[(z - M_k)^{-2}, \Omega_{N,k}] E_k[(\bar{w} - M_k^*)^{-2}, \Omega_{N,k}]) e_k,$$

in which we have temporarily removed the notation and written everything out in full for the sake of clarity.

The next step is to replace the question of the convergence of (3.11) with that of its antiderivative. In particular, what we proceed to prove is that

(3.12)
$$\frac{1}{N} \sum_{k=1}^{N} e_k^T (E_k[\Xi_{N,k}(z)^{-1}, \Omega_{N,k}] E_k[(\Xi_{N,k}(w)^{-1})^*, \Omega_{N,k}]) e_k$$
$$\to \log\left(1 - \frac{1}{z\bar{w}}\right)$$

in probability. That this is sufficient follows from a theorem of Vitali (see, e.g., [29]): If $\{f_N(z)\}$ is a sequence of uniformly bounded analytic functions on a domain $\mathbb{D} \subset C$, converging at a collection $\{z\}$ containing a limit point in $\mathbb{D}$, then there exists an analytic $f(z)$ such that $f_N(z) \to f(z)$ and $f'_N(z) \to f'(z)$ throughout $\mathbb{D}$. Denote the left-hand side of (3.12) by $\phi_N(z, \bar{w})$. This function is analytic in $z$ ($|z| > 4$) for any fixed $\bar{w}$ ($|\bar{w}| > 4$); in the same region, it is analytic in $\bar{w}$ for any fixed $z$. That $\phi_N(z, \bar{w})$ is bounded independently of $N$ follows from the fact that

(3.13) $\|E_k[\Xi_{N,k}(z)^{-1}, \Omega_{N,k}]\|^2 \leq E_k[\|\Xi_{N,k}(z)^{-1}\|^2, \Omega_{N,k}] \leq \mathcal{K}^2.$

Now, if the convergence (3.12) takes place for given $z$ and $\bar{w}$, we can fix a collection $\{z_k, \bar{w}_k\}$ with limit point in $\{|z|, |\bar{w}| > 4\}$ and a subsequence $N' \to \infty$ over which the convergence is almost sure, simultaneously, for all points $(z_k, \bar{w})$. On such a realization, the theorem of Vitali can be applied twice, one coordinate at a time, to obtain the convergence of $\phi_N$ and its mixed partials over the full domain. Since the conclusion is independent of both the subsequence and the particular realization, (3.12) does indeed imply the latter half of (3.9).

Moving on, we bring in some additional shorthand,

$$\mathcal{A}_{N,k}(z) \equiv \Xi_{N,k}(z)^{-1}\mathbf{1}_{\Omega_{N,k}}, \qquad \mathcal{A}_{N,k\ell}(z) \equiv \Xi_{N,k\ell}(z)^{-1}\mathbf{1}_{\Omega_{N,k\ell}}, \qquad \text{and so on,}$$

and focus on a fixed term of (3.11),

$$(3.14) \qquad T_{N,k}(z,\bar{w}) = e_k^T E_k[\mathcal{A}_{N,k}(z)] E_k[\mathcal{A}_{N,k}(w)^*] e_k.$$

The method is rather laborious. First, successive columns of $M$ are removed [effectively replacing $\mathcal{A}_{N,k}(z)$ with $\mathcal{A}_{N,k\ell}(z)$, etc.] in order to invoke the basic estimate of Lemma 2.3. Afterward, the removed columns are effectively restored, arriving at an equation for $T_{N,k}(z,\bar{w})$ which is approximately deterministic for large $N$.

Invoking the matrix identity

$$(3.15) \qquad \Xi_{N,k}(z)^{-1} = \frac{1}{z}I + \frac{1}{z}\sum_{\ell \neq k} \Xi_{N,k}(z)^{-1} m_\ell e_\ell^T,$$

(3.14) is expanded: with $P_k$ denoting the conditional probability with respect to $m_1, \ldots, m_k$, we have

$$
\begin{aligned}
T_{N,k}(z,\bar{w}) \\
&= \frac{1}{z\bar{w}} P_k^2(\Omega_{N,k}) + \frac{1}{z\bar{w}} \sum_{\ell \neq k} (e_k^T E_k[\mathcal{A}_{N,k}(z) m_\ell] e_\ell^T e_k) P_k(\Omega_{N,k}) \\
&\quad + \frac{1}{z\bar{w}} \sum_{\ell \neq k} e_k^T e_\ell E_k[m_\ell^* \mathcal{A}_{N,k}(w)^*] e_k \\
&\quad + \frac{1}{z\bar{w}} \sum_{\ell \neq k} \sum_{j \neq k} e_k^T E_k[\mathcal{A}_{N,k}(z) m_\ell] e_\ell^T e_j E_k[m_j^* \mathcal{A}_{N,k}(w)^*] e_k.
\end{aligned}
$$
(3.16)

This, in turn, may immediately be reduced to

$$
\begin{aligned}
T_{N,k}(z,\bar{w}) &= \frac{1}{z\bar{w}} + \frac{1}{z\bar{w}} \sum_{\ell \neq k} e_k^T E_k[\mathcal{A}_{N,k}(z) m_\ell] E_k[m_\ell^* \mathcal{A}_{N,k}(w)^*] e_k \\
&\quad + \mathcal{E}_{N,k}^0(z,\bar{w}).
\end{aligned}
$$
(3.17)

Here, $\mathcal{E}_{N,k}^0(z,\bar{w}) = (z\bar{w})^{-1}(1 - P_k^2(\Omega_{N,k}))$, which is easily seen to go to zero in $L^1$. Note also that the second and third terms of (3.16) are identically zero: both contain $e_k^T e_\ell = 0$ in each $\ell \neq k$ term of their defining sums. The reduction of the last term of (3.16) to the sum present in (3.17) follows similarly. As for that sum, we wish to replace the $\mathcal{A}_{N,k}(\cdot)$'s with the corresponding $\mathcal{A}_{N,k\ell}(\cdot)$. The procedure and goal, which is to produce independence between the resolvent matrices and the columns of $M$, follows the outline of the proof of Lemma 3.1.

First, note that for $\ell \neq k$,

$$\Xi_{N,k}(z)^{-1} m_\ell = \delta_{k\ell}(z) \Xi_{N,k\ell}(z)^{-1} m_\ell,$$

in which

$$\delta_{k\ell}(z) \equiv (1 + e_\ell^T \Xi_{N,k\ell}(z)^{-1} m_\ell)^{-1} \equiv (1 + \varepsilon_{k\ell}(z))^{-1}$$

[recall (3.3)]. Now write

$$\sum_{\ell \neq k} e_k^T E_k[\mathcal{A}_{N,k}(z) m_\ell] E_k[m_\ell^* \mathcal{A}_{N,k}(w)^*] e_k$$

$$= \sum_{\ell \neq k} e_k^T E_k[\Xi_{N,k\ell}(z)^{-1} m_\ell, \Omega_{N,k}]$$

(3.18)

$$\times E_k[m_\ell^*(\Xi_{N,k\ell}(w)^{-1})^*, \Omega_{N,k}] e_k + \mathcal{E}_{N,k}^1$$

$$= \sum_{\ell < k} e_k^T E_k[\mathcal{A}_{N,k\ell}(z)](m_\ell m_\ell^*) E_k[\mathcal{A}_{N,k\ell}(w)^*] e_k + \mathcal{E}_{N,k}^1 + \mathcal{E}_{N,k}^2,$$

where $E_k[\mathcal{A}_{N,k,\ell}(\cdot) m_\ell] = 0$ for $\ell > k$ is used in the last line. $\mathcal{E}_{N,k}^1 = \mathcal{E}_{N,k}^1(z, \bar{w})$ and $\mathcal{E}_{N,k}^2 = \mathcal{E}_{N,k}^2(z, \bar{w})$ are the sums over $\ell \neq k$ of

$$\mathcal{E}_{N,k,\ell}^1(z, \bar{w}) = E_k[(\delta_{k\ell}(z) - 1) e_k^T \Xi_{N,k\ell}(z)^{-1} m_\ell, \Omega_{N,k}]$$

$$\times E_k[m_\ell^*(\Xi_{N,k\ell}(w)^{-1})^* e_k, \Omega_{N,k}]$$

$$+ E_k[e_k^T \Xi_{N,k\ell}(z)^{-1} m_\ell, \Omega_{N,k}]$$

(3.19)

$$\times E_k[m_\ell^*(\Xi_{N,k\ell}(w)^{-1})^* e_k (\overline{\delta_{k\ell}(w)} - 1), \Omega_{N,k}]$$

$$+ E_k[(\delta_{k\ell}(z) - 1) e_k^T \Xi_{N,k\ell}(z)^{-1} m_\ell, \Omega_{N,k}]$$

$$\times E_k[m_\ell^*(\Xi_{N,k\ell}(w)^{-1})^* e_k (\overline{\delta_{k\ell}(w)} - 1) \Omega_{N,k}]$$

and

$$\mathcal{E}_{N,k,\ell}^2(z, \bar{w}) = -E_k[e_k^T \mathcal{A}_{N,k\ell}(z) m_\ell] E_k[m_\ell^* \mathcal{A}_{N,k\ell}(w)^* e_k, \Omega_{N,k}^C]$$

(3.20)

$$- E_k[e_k^T \mathcal{A}_{N,k\ell}(z) m_\ell, \Omega_{N,k}^C] E_k[m_\ell^* \mathcal{A}_{N,k\ell}(w)^* e_k]$$

$$+ E_k[e_k^T \mathcal{A}_{N,k\ell}(z) m_\ell, \Omega_{N,k}^C] E_k[m_\ell^* \mathcal{A}_{N,k\ell}(w)^* e_k, \Omega_{N,k}^C].$$

Similarly to (3.7), that is, by Lemma 2.3, we have

(3.21)
$$E|E_k[e_k^T \Xi_{N,k\ell}(\cdot)^{-1} m_\ell, \Omega_{N,k}]|^2 \leq E[|e_k^T \Xi_{N,k\ell}(\cdot)^{-1} m_\ell|^2, \Omega_{N,k\ell}]$$
$$= E|m_\ell^* \mathcal{A}_{N,k,\ell}(\cdot)^* e_k e_k^T \mathcal{A}_{N,k,\ell}(\cdot) m_\ell|$$
$$\leq C_1 N^{-1}$$

for $\cdot = z, w$. Next, note the equality

$$(\delta_{k\ell}(\cdot) - 1)(e_k^T \Xi_{N,k\ell}(\cdot)^{-1} m_\ell) = (-\varepsilon_{k\ell}(\cdot) + \varepsilon_{k\ell}(\cdot)^2 \delta_{k\ell}(\cdot))(e_k^T \Xi_{N,k\ell}(\cdot)^{-1} m_\ell)$$
$$= -(e_\ell^T \Xi_{N,k\ell}(\cdot)^{-1} m_\ell)(e_k^T \Xi_{N,k\ell}(\cdot)^{-1} m_\ell)$$
$$+ (e_\ell^T \Xi_{N,k\ell}(\cdot)^{-1} m_\ell)^2 (e_k^T \Xi_{N,k}(\cdot)^{-1} m_\ell),$$

the identity $\delta_{k\ell}(\cdot)\Xi_{N,k\ell}(\cdot)^{-1} m_\ell = \Xi_{N,k}(\cdot)^{-1} m_\ell$ having been reinvoked in the last step. It then follows by the same type of reasoning used in (3.21) that

(3.22) $\quad E|E_k[(\delta_{k\ell}(\cdot) - 1) e_k^T \Xi_{N,k\ell}(\cdot)^{-1} m_\ell, \Omega_{N,k}]|^2 \leq 2\mathcal{I}_1 + 2\mathcal{I}_2,$

where

$$
\begin{aligned}
(3.23) \quad \mathcal{I}_1 &= E[|e_\ell^T \Xi_{N,k\ell}(\cdot)^{-1} m_\ell|^2 |e_k^T \Xi_{N,k\ell}(\cdot)^{-1} m_\ell|^2, \Omega_{N,k}] \\
&\leq E[|e_\ell^T \mathcal{A}_{N,k\ell}(\cdot) m_\ell|^2 |e_k^T \mathcal{A}_{N,k\ell}(\cdot) m_\ell|^2] \\
&\leq C_2 N^{-2},
\end{aligned}
$$

after an application of the Schwarz inequality and

$$
\begin{aligned}
(3.24) \quad \mathcal{I}_2 &= E[|e_\ell^T \Xi_{N,k\ell}(\cdot)^{-1} m_\ell|^4 |e_k^T \Xi_{N,k}(\cdot)^{-1} m_\ell|^2, \Omega_{N,k}] \\
&\leq E[|e_\ell^T \Xi_{N,k\ell}(\cdot)^{-1} m_\ell|^4 \|\Xi_{N,k}(\cdot)^{-1} e_k\|_2^2 \|m_\ell\|_2^2, \Omega_{N,k}] \\
&\leq \mathcal{K}^2 E[|e_\ell^T \mathcal{A}_{N,k\ell}(\cdot) m_\ell|^4 \|m_\ell\|_2^2] \\
&\leq \mathcal{K}^2 (E[\|m_\ell\|_2^4])^{1/2} E[|e_\ell^T \mathcal{A}_{N,k\ell}(\cdot) m_\ell|^8]^{1/2} \\
&\leq C_3 N^{-2}.
\end{aligned}
$$

Here, $\|\cdot\|_2$ denotes the usual $\ell^2$-norm.

Combining (3.21) and (3.22)–(3.24) gives

$$E|\mathcal{E}_{N,k}^1(z, \bar{w})| = O(N^{-1/2}).$$

On the other hand, (3.21) plus the sharp decay $P(\Omega_{N,k}^C)$ is enough to produce a similar bound on the $L^1$-norm of $\mathcal{E}_{N,k}^2(z, \bar{w})$.

Returning to (3.18), since $E[m_\ell m_\ell^*] = (1/N)I$, it is natural to claim that

$$
(3.25) \quad \begin{aligned}
&\sum_{\ell < k} e_k^T E_k[\mathcal{A}_{N,k\ell}(z)](m_\ell m_\ell^*) E_k[\mathcal{A}_{N,k\ell}(w)^*] e_k \\
&\qquad = \frac{1}{N} \sum_{\ell < k} e_k^T E_k[\mathcal{A}_{N,k\ell}(z)] E_k[\mathcal{A}_{N,k\ell}(w)^*] e_k + \mathcal{E}_{N,k}^3(z, \bar{w})
\end{aligned}
$$

with a small error; we show below that

(3.26) $\qquad\qquad\qquad E|\mathcal{E}_{N,k}^3(z, \bar{w})| = o(1).$

At this point, the endgame becomes clear. If we could now restore the last column removed from the resolvent through an estimate along the lines of

$$
\begin{aligned}
(3.27) \quad E|e_k^T(E_k[\mathcal{A}_{N,k\ell}(z)]E_k[\mathcal{A}_{N,k\ell}(w)^*] \\
- E_k[\mathcal{A}_{N,k}(z)]E_k[\mathcal{A}_{N,k}(w)^*])e_k| = o(1),
\end{aligned}
$$

the combined outcome of this statement with (3.17), (3.18) and (3.25) would be that

$$ T_{N,k}(z,\bar{w}) = \frac{1}{z\bar{w}} + \frac{1}{z\bar{w}}\left(\frac{k-1}{N}T_{N,k}(z,\bar{w})\right) + \mathcal{E}_{N,k}(z,\bar{w}), $$

in which

$$ E|\mathcal{E}_{N,k}(z,\bar{w})| = o(1), $$

uniformly in $k$, $1 \le k \le N$. This is to say that the object under investigation, the left-hand side of (3.12), satisfies

$$ \left|\frac{1}{N}\sum_{k=1}^N T_{N,k}(z,\bar{w}) - \frac{1}{N}\sum_{k=1}^N \left(z\bar{w} - \frac{k-1}{N}\right)^{-1}\right| \le \left(\frac{|z\bar{w}|}{|z\bar{w}|-1}\right)\frac{1}{N}\sum_{k=1}^N |\mathcal{E}_{N,k}(z,\bar{w})| $$

and so tends to zero in probability. The advertised limit follows from identifying the Riemann sum.

Returning to the bounds (3.26) and (3.27), we begin with the second. The arguments used in controlling the error $\mathcal{E}_{N,k}^2(z,\bar{w})$ above will show why

$$
\begin{aligned}
E|e_k^T E_k[\mathcal{A}_{N,k\ell}(z)]E_k[\mathcal{A}_{N,k\ell}(w)^*]e_k \\
- e_k^T E_k[\Xi_{N,k\ell}(z)^{-1},\Omega_{N,k}]E_k[(\Xi_{N,k\ell}(w)^{-1})^*,\Omega_{N,k}]e_k|
\end{aligned}
$$

tends to zero faster than, say, $N^{-1}$. Next, to restore $\Xi_{N,k\ell}(\cdot)^{-1}$ to $\Xi_{N,k}(\cdot)^{-1}$, note, for example, that

$$
\begin{aligned}
E|E_k[e_k^T(\Xi_{N,k}(z)^{-1} - \Xi_{N,k\ell}(z)^{-1}),\Omega_{N,k}]E_k[(\Xi_{N,k}(w)^{-1})^*e_k,\Omega_{N,k}]|^2 \\
\le E[E_k[|e_k^T\Xi_{N,k\ell}(z)^{-1}m_\ell|^2\|(\Xi_{N,k}(z)^{-1})^*e_\ell\|_2^2,\Omega_{N,k}] \\
\times E_k[\|(\Xi_{N,k}(w)^{-1})^*e_k\|_2^2,\Omega_{N,k}]] \\
\le \mathcal{K}^4 E|e_k^T \mathcal{A}_{N,k\ell}(z)m_\ell|^2 \\
= O(N^{-1}).
\end{aligned}
$$

It will follow that the right-hand side of (3.27) is of order $N^{-1/2}$ for any $k, \ell$.

Finally, to verify (3.26), we take the mean-square of that error term, written as

$$ \frac{1}{N^2}E\left|\sum_{\ell<k}(\tilde{m}_\ell^* E_k[\mathcal{A}_{N,k\ell}(w)^*](e_k e_k^T)E_k[\mathcal{A}_{N,k\ell}(z)]\tilde{m}_\ell \right. $$

(3.28)
$$-e_k^T E_k[\mathcal{A}_{N,k\ell}(z)] E_k[\mathcal{A}_{N,k\ell}(w)^*] e_k) \bigg|^2$$

(3.29)
$$\equiv \frac{1}{N^2} E \bigg| \sum_{\ell < k} (\tilde{m}_\ell^* D_{N,k\ell}(z,\bar{w}) \tilde{m}_\ell - \operatorname{tr}(D_{N,k\ell}(z,\bar{w}))) \bigg|^2,$$

where $\tilde{m}_\ell = \sqrt{N} m_\ell$ are the columns of the unnormalized matrix. Again, by Lemma 2.3, we have that

(3.30)
$$E |\tilde{m}_\ell^* D_{N,k\ell}(z,\bar{w}) \tilde{m}_\ell - \operatorname{tr}(D_{N,k\ell}(z,\bar{w}))|^2 \leq C_4$$

for a fixed constant $C_4$. Therefore, the contribution from the diagonal terms of (3.28) is of order $N^{-1}$. Consideration of the off-diagonal terms is complicated by the dependence between $E_k[\mathcal{A}_{N,k\ell}(\cdot)]$ and $E_k[\mathcal{A}_{N,kj}(\cdot)]$, forcing us to remove yet another column to arrive at the required estimate. Note that with $j, \ell < k$, the typical cross term satisfies

$$E[(\tilde{m}_\ell^* D_{N,k\ell}(z,\bar{w}) \tilde{m}_\ell - \operatorname{tr}(D_{N,k\ell}(z,\bar{w})))$$
$$\times (\tilde{m}_j^* D_{N,kj}(z,\bar{w}) \tilde{m}_j - \operatorname{tr}(D_{N,kj}(z,\bar{w})))]$$
$$= E[(\tilde{m}_\ell^* D_{N,k\ell}(z,\bar{w}) \tilde{m}_\ell - \operatorname{tr}(D_{N,k\ell}(z,\bar{w})))$$
$$\times (\tilde{m}_j^* (D_{N,kj}(z,\bar{w}) - D_{N,kj\ell}(z,\bar{w})) \tilde{m}_j$$
$$- \operatorname{tr}(D_{N,kj}(z,\bar{w}) - D_{N,kj\ell}(z,\bar{w})))],$$

where in $D_{N,kj\ell}(\cdot,\cdot)$, each $\mathcal{A}_{N,kj}(\cdot)$ is replaced by $\mathcal{A}_{N,kj\ell}(\cdot)$. After two applications of the Schwarz inequality and (3.30), the required estimate will follow if

(3.31)
$$E|\tilde{m}_j^* (D_{N,kj}(z,\bar{w}) - D_{N,kj\ell}(z,\bar{w})) \tilde{m}_j|^2 = o(1)$$

and

(3.32)
$$E|\operatorname{tr} D_{N,kj}(z,\bar{w}) - \operatorname{tr} D_{N,kj\ell}(z,\bar{w})|^2 = o(1)$$

for $N \to \infty$. However, Lemma 2.3 bounds the left-hand side of (3.31) in terms of the left-hand side of (3.32). It remains to observe that after unraveling the notation, (3.32) is really the same as (3.27). We need to compare combinations of the resolvent with two columns removed, $e_k^T E_k[\mathcal{A}_{N,k\ell}(z)] E_k[\mathcal{A}_{N,k\ell}(w)^*] e_k$, to the same object in which a third column has been removed. (3.27) compares that same combination with one column removed to its image after a second column has been removed. An identical procedure used to bound (3.27) will then show that (3.32) is $O(N^{-1/2})$. The proof is thus complete. □

3.2. *Tightness estimate.* The previous section established the convergence in distribution of any finite collection $\{\mathfrak{G}_N^0(z_1), \mathfrak{G}_N^0(z_2), \ldots, \mathfrak{G}_N^0(z_k)\}$ of the truncated and recentered resolvent

$$\mathfrak{G}_N^0(\cdot) = \operatorname{tr} \Xi_N(\cdot)^{-1} \mathbf{1}_{\Omega_N} - E[\operatorname{tr} \Xi_N(\cdot)^{-1}, \Omega_N]$$
$$= \sum_{k=1}^N (E_k - E_{k-1})[\operatorname{tr} \Xi_N(z)^{-1}, \Omega_N],$$

sampled at any points $z_\ell, 1 \leq \ell \leq k$, exterior to the disk of radius four. The standard method for extending such finite-dimensional convergence to the convergence of $\mathfrak{G}_N^0(z)$ in the space of continuous processes over a predetermined contour $\mathcal{C}$ is to check the Arzela–Ascoli criteria (see, e.g., [5]): verifying the estimate

$$E|\mathfrak{G}_N^0(z) - \mathfrak{G}_N^0(w)|^\alpha \leq C|z - w|^{1+\beta}$$

for any $\alpha, \beta > 0$ and all $z, w \in \mathcal{C}$ is enough.

In the present case, it is convenient to recenter once more before preceding. Note that

$$\sum_{k=1}^N (E_k - E_{k-1})[\operatorname{tr} \Xi_{N,k}(z)^{-1}] = 0$$

for any $z$, while for any $\mathcal{C}$ exterior to the disk of radius four,

$$\sup_{z \in \mathcal{C}} \sum_{k=1}^N (E_k - E_{k-1})[\operatorname{tr} \Xi_{N,k}(z)^{-1}, \Omega_N^C] \to 0$$

in probability. The proof of the latter is the same as that for $\sup_{z \in \mathcal{C}} |\mathfrak{G}_N^1(z)| \to 0$ in Section 2. It follows that the tightness of $\mathfrak{G}_N^0(\cdot)$ in the space of continuous functions on a given contour $\mathcal{C}$ is equivalent to that of

$$\widetilde{\mathfrak{G}}_N^0(z) \equiv \sum_{k=1}^N (E_k - E_{k-1})[\operatorname{tr}(\Xi_N(z)^{-1} - \Xi_{N,k}(z)^{-1}), \Omega_N].$$

The following estimate, along with Lemma 4.1 proved in the next section, completes the proof of Theorem 1.2:

LEMMA 3.3. *It holds that*

$$E\left|\frac{\widetilde{\mathfrak{G}}_N^0(z) - \widetilde{\mathfrak{G}}_N^0(z)}{z - w}\right|^2 \leq C,$$

*with a constant $C$ independent of both $N$ and any choice of $z, w$ such that $|z|, |w| > 4$.*

PROOF. It is no surprise that considerations are much the same as those behind the proofs of Lemmas 3.1 and 3.2. We at once have that

$$(w-z)^{-1}(\widetilde{\mathfrak{G}}_N^0(z) - \widetilde{\mathfrak{G}}_N^0(w))$$
$$= \sum_{k=1}^N (E_k - E_{k-1})$$
$$\times [\mathrm{tr}(\Xi_N(z)^{-1}\Xi_N(w)^{-1} - \Xi_{N,k}(z)^{-1}\Xi_{N,k}(w)^{-1}), \Omega_N].$$

Now, again with $\delta_k(\cdot) = (1 + e_k^T \Xi_{N,k}(\cdot)^{-1} m_k)^{-1}$, we write

$$\Xi_N(z)^{-1}\Xi_N(w)^{-1} - \Xi_{N,k}(z)^{-1}\Xi_{N,k}(w)^{-1}$$
$$= (\Xi_N(z)^{-1} - \Xi_{N,k}(z)^{-1})(\Xi_N(w)^{-1} - \Xi_{N,k}(w)^{-1})$$
$$+ (\Xi_N(z)^{-1} - \Xi_{N,k}(z)^{-1})\Xi_{N,k}(w)^{-1}$$
$$+ \Xi_{N,k}(z)^{-1}(\Xi_N(w)^{-1} - \Xi_{N,k}(w)^{-1})$$
$$= \delta_k(z)\delta_k(w)(\Xi_{N,k}(z)^{-1}(m_k e_k^T)$$
$$\times \Xi_{N,k}(z)^{-1}\Xi_{N,k}(w)^{-1}(m_k e_k^T)\Xi_N(w)^{-1})$$
$$+ \delta_k(z)(\Xi_{N,k}(z)^{-1}(m_k e_k^T)\Xi_{N,k}(z)^{-1}\Xi_{N,k}(w)^{-1})$$
$$+ \delta_k(w)(\Xi_{N,k}(z)^{-1}\Xi_{N,k}(w)^{-1}(m_k e_k^T)\Xi_{N,k}(w)^{-1})$$

and so

$$\mathrm{tr}(\Xi_N(z)^{-1}\Xi_N(w)^{-1} - \Xi_{N,k}(z)^{-1}\Xi_{N,k}(w)^{-1})$$
$$= \delta_k(z)\delta_k(w)(e_k^T \Xi_{N,k}(z)^{-1}\Xi_k(w)^{-1} m_k)^2$$
(3.33)
$$+ \delta_k(z)(e_k^T \Xi_{N,k}(z)^{-2}\Xi_{N,k}(w)^{-1} m_k)$$
$$+ \delta_k(w)(e_k^T \Xi_{N,k}(w)^{-2}\Xi_{N,k}(z)^{-1} m_k).$$

We therefore wish to show that

(3.34) $\quad E|\delta_k(z)\delta_k(w)(e_k^T \Xi_{N,k}(z)^{-1}\Xi_{N,k}(w)^{-1} m_k)^2, \Omega_N|^2 \leq C_1 N^{-1}$

and

(3.35) $\quad E|\delta_k(z)(e_k^T \Xi_{N,k}(z)^{-2}\Xi_{N,k}(w)^{-1} m_k), \Omega_N|^2 \leq C_1 N^{-1}$

for arbitrary $k$, $z$ and $w$ [the last two terms of (3.33) are the same up to exchanging $z$ and $w$].

We provide the full details behind (3.34) only, this being the more elaborate and both assertions requiring the same strategy. Further, the strategy is familiar. In particular, the reader will recognize that the desired bounds follow directly from Lemma 2.3 if $\delta_k(z)$ and $\delta_k(w)$ can be replaced with 1 inside both expectations:

$$E|(e_k^T \Xi_{N,k}(z)^{-1}\Xi_{N,k}(w)^{-1}m_k)^2, \Omega_N|^2$$
$$\leq E[|m_k^*(\Xi_{N,k}(z)^{-1}\Xi_{N,k}(w)^{-1})^*(e_k e_k^T)(\Xi_{N,k}(z)^{-1}\ \Xi_{N,k}(w)^{-1})m_k|^2,$$
$$\Omega_{N,k}]$$
(3.36)
$$\leq C_2 N^{-2} E[|e_k^T\ (\Xi_{N,k}(z)^{-1}\Xi_{N,k}(w)^{-1})(\Xi_{N,k}(z)^{-1}\Xi_{N,k}(w)^{-1})^* e_k|^2,$$
$$\Omega_{N,k}]$$
$$\leq C_2 \mathcal{K}^8 N^{-2}$$

and, likewise,

$$E|e_k^T \Xi_{N,k}(z)^{-2}\Xi_{N,k}(w)^{-1}m_k, \Omega_N|^2$$
$$\leq E[|m_k^*(\Xi_{N,k}(z)^{-2}\Xi_{N,k}(w)^{-1})^*(e_k e_k^T)(\Xi_{N,k}(z)^{-2}\Xi_{N,k}(w)^{-1})m_k|,$$
$$\Omega_{N,k}]$$
$$\leq C_3 \mathcal{K}^6 N^{-1}.$$

Recall the notation $\delta_k(z) = (1 + \varepsilon_k(z))^{-1}$ and set

$$\delta_k(z)\delta_k(w) = (1 + \varepsilon_k(z) + \varepsilon_k(w) + \varepsilon_k(z)\varepsilon_k(w))^{-1} = (1 + \gamma_k(z,w))^{-1}.$$

The key step in proving (3.34) is to notice the a priori bound: on the event $\Omega_N$,

(3.37)
$$|\delta_k(z)\delta_k(w)(e_k^T \Xi_{N,k}(z)^{-1}\Xi_{N,k}(w)^{-1}m_k)^2|$$
$$= |\mathrm{tr}(\Xi_N(z)^{-1} - \Xi_{N,k}(z)^{-1})(\Xi_N(w)^{-1} - \Xi_{N,k}(w)^{-1})|$$
$$\leq N(\|\Xi_N(z)^{-1}\| + \|\Xi_{N,k}(z)^{-1}\|)(\|\Xi_N(w)^{-1}\| + \|\Xi_{N,k}(w)^{-1}\|)$$
$$\leq 4\mathcal{K}^2 N.$$

Next, we know by now that $E[|\varepsilon_k(\cdot)|^{2p}, \Omega_{N,k}] \leq C_4 N^{-p}$ for $p \geq 1$ and thus there is a constant $C_5$ such that

(3.38) $$E[|\gamma_k(z,w)|^{2p}, \Omega_N] \leq C_5 N^{-p}.$$

If we now write

$$\delta_k(z)\delta_k(w) = 1 - \gamma_k(z,w) + \gamma_k(z,w)^2 - \gamma_k(z,w)^3 \delta_k(z)\delta_k(w)$$

in the left-hand side of (3.34), successive applications of the Schwarz inequality, along with (3.36) and (3.38), will show that each of the resulting first three terms are of order $N^{-2}$. For the fourth term, the bound (3.37), followed by the appraisal $E[|\gamma_k(z,w)|^6, \Omega_N] = O(N^{-3})$, will produce the order $N^{-1}$ in (3.34). For (3.35), the argument is as follows. Playing a similar role to (3.37), the fact that $\delta_k(z)(\Xi_{N,k}(z)^{-2}\Xi_k(w)^{-1}m_k)\mathbf{1}_{\Omega_N}$ is less than a multiple of $N$ is noted. Then $\delta_k(z)$ is itself expanded to third order with each term being handled in an identical manner to that just described for (3.35). The proof is thus complete. □

**4. Equivalence of the means.** We prove the following:

LEMMA 4.1. *It holds that*
$$\lim_{N\to\infty} |E[\operatorname{tr} \Xi_N(z)^{-1}] - Nz^{-1}| = 0,$$
*uniformly for $z$ belonging to any contour exterior to the unit disk.*

In the notation of Section 2, this reads as $|\mathfrak{G}_N^2(z)| \to 0$. Therefore, the central limit theorem is valid for linear statistics centered by the asymptotic mean as well as the $N$th level mean. Note that the above result requires only that $|z| > 1$. As we have discussed, this is certainly the appropriate (i.e., optimal) condition for Theorem 1.2. While the technique used to prove the CLT rests on comparing rank one perturbations of the random matrix $M$, the present computation requires no such comparisons. Thus, rather than relying on bounds on the spectral norm to provide a priori control, we are able to use the spectral radius which is more natural to the problem.

PROOF OF LEMMA 4.1. Denote
$$\Lambda_N = \max\{|\lambda| : \lambda \text{ is an eigenvalue of } M\}.$$

Along with (2.1), Geman [13] has also proved that for any $\rho > 1$ and $\alpha > 0$,

(4.1) $$P(\Lambda_N > \rho) = o(N^{-\alpha}).$$

(Once again, [13] discusses the case that the entries of the matrix are real, but the extension to complex entries is effortless.) Also from [13] [see the display after equation (6) of that paper], we record the bound

(4.2) $$E[\Lambda_N^p] \leq C_1 N^2 p$$

for a constant $C_1$ and $p$ allowed to grow with $N$, up to $p = O(\ln N)$.

Now, fix a positive $\rho$ larger than one but less than $|z|$. The bound (4.1), together with the often-used Lemma 2.2, implies that
$$E[\operatorname{tr} \Xi_N(z)^{-1}, \Lambda_N > \rho] \to 0$$
as $N \to \infty$. On the opposite event, $\Xi_N(z)^{-1}$ may be expanded inside the expectation as in
$$E[\Xi_N(z)^{-1}, \Lambda_N \leq \rho]$$
$$= \sum_{k=0}^{p-1} \frac{1}{z^{k+1}} E[\operatorname{tr} M^k, \Lambda_N \leq \rho] + \frac{1}{z^{p+1}} E[\operatorname{tr}(\Xi_N(z)^{-1} M^p), \Lambda_N \leq \rho].$$

The first term in the sum on the right-hand side reads $z^{-1}(N - NP(\lambda_N < \rho)) = Nz^{-1} + o(1)$, in which we see the desired result. Next, note that the remainder term is bounded by $N(|z| - \rho)^{-1}(\rho/|z|)^p$ and thus goes to zero for $p = p(N) = C_2 \ln N$ with a sufficiently large constant $C_2$. Further, on account of (4.2),

$$E[\operatorname{tr} M^k, \Lambda_N > \rho] \leq NE[\Lambda_N^{2k}]^{1/2} P(\Lambda_N > \rho)^{1/2} \leq C_3 N^{-1}$$

for all $N$ sufficiently large and all $k < (C_2/2) \ln N$. Therefore, it remains to prove an estimate of the type

$$(4.3) \qquad \lim_{N \to \infty} \sum_{1 \leq s \leq C_1 \ln N} E[\operatorname{tr} M^s] = 0.$$

By the assumptions that $E[m_{k\ell}] = E[m_{k\ell}^2] = 0$, we have

$$|E \operatorname{tr} M^s| \leq \sum\nolimits^* E|m_{\ell_1 \ell_2} m_{\ell_2 \ell_3} \cdots m_{\ell_s \ell_1}|,$$

in which $\sum^*$ denotes the sum over those choices of indices $\{\ell_1, \ell_2, \ldots, \ell_s\}$ such that any index pair $(\ell_j, \ell_{j+1})$, or *entry type* $m_{\ell_j \ell_{j+1}}$, which appears must do so at least three times. This object may be bounded in a crude but sufficient manner as follows:

$$\sum\nolimits^* E|m_{\ell_1 \ell_2} m_{\ell_2 \ell_3} \cdots m_{\ell_s \ell_1}|$$
$$\leq \sum_{k=1}^{\lfloor s/3 \rfloor} \sum_{n_1 + \cdots + n_k = s} \frac{s!}{n_1! n_2! \cdots n_k!} N^k (E|m_{11}|^{n_1} E|m_{11}|^{n_2} \cdots E|m_{11}|^{n_k}).$$

The outermost sum indicates that $k = 1, \ldots, \lfloor s/3 \rfloor$ distinct entry types can be present in any string being summed over. The second sum is over all possible numbers of each entry type: $n_1$ type 1 entries, $n_2$ type 2 entries, and so on. We allow for all nonnegative, ordered $k$-tuples $n_1, n_2, \ldots, n_k$ summing to $s$, which results in an overestimate. The next factor, $s!/n_1! \cdots n_k!$, assigns each entry from left to right one of the $k$ types. Afterward, each type must be labeled by choosing the corresponding pair of indices. That is, the first entry in the string might have been assigned "type 3" at this point, but type 3 can be $m_{1,2}$, or $m_{7,4}$, and so on. Note that not all such unrestricted type assignments will support a consistent labeling of indices; we have another overestimate. As for the number of ways to label the indices, the factor $N^k$ is an obvious upper bound. Moving from left to right, each time a new (unassigned) entry type is encountered, one of $N$ possible choices for its first index is chosen. After the first index corresponding to each particular type is chosen, all of the second indices are determined by their neighboring entries. Again, we allow for the overcount produced by a number of inconsistent labelings which result from this procedure. The final product of expectations

requires no explanation and by our assumptions [see (1.2.iii)], it is bounded as in

$$E|m_{11}|^{n_1} E|m_{11}|^{n_2} \cdots E|m_{11}|^{n_k} \leq N^{-s/2} \prod_{\ell=1}^{k} (n_\ell)^{\alpha n_\ell} \leq N^{-s/2} s^{\alpha s}$$

for some $\alpha > 0$. It follows that

$$|E \operatorname{tr} M^s| \leq N^{-s/2} s^{\alpha s} \sum_{k=1}^{\lfloor s/3 \rfloor} N^k k^s \leq N^{-s/6} s^{((\alpha+1)s+1)}.$$

This clearly goes to zero as $N \to \infty$ with $s = O(\ln N)$, and it does so fast enough that the sum over all such $s$ will also vanish [recall (4.3)]. The proof is thus complete. $\square$

## APPENDIX

We include here a proof of Lemma 2.2. Important use is made of the following fact which was established in [1] for a similar purpose:

PROPOSITION A.1 (Corollary A.2 of [1]). *Let the vector $X = (x_1, x_2, \ldots, x_n)$ consist of independent complex entries such that the joint density of the real and imaginary parts of each $x_\ell$ are uniformly bounded by a constant $C$. Let $\alpha_1, \alpha_2, \ldots, \alpha_k$ be $k$ orthonormal complex unit vectors. Then the joint density of the complex numbers $X^*\alpha_1, X^*\alpha_2, \ldots, X^*\alpha_k$ is bounded by $C^{2k} n^{2k}$.*

PROOF OF LEMMA 2.2. For any $N \times N$ matrix $A$, $|\sum_{k=1}^{N} \lambda_k(A)|$ is dominated by $\sum_{k=1}^{N} [(A^*A)_{kk}]^{1/2}$ and so

$$|\operatorname{tr} \Xi_N(z)^{-1}|^p \leq N^{p-1} \sum_{k=1}^{N} [(\Xi_N(z)^* \Xi_N(z))^{-1}_{kk}]^{p/2}$$

when $p > 1$. The notation on the right indicates (the $p/2$ power of) the $kk$ entry of the inverse matrix of $\Xi_N(z)^* \Xi_N(z)$. Next, let the $k$th column of $\Xi_N(z)$ be denoted by $\xi_k(z)$. Also, let $\Xi_{N(k)}(z)$ denote the $N \times (N-1)$ matrix formed by removing $\xi_k(z)$ from $\Xi_N(z)$. Then

(A.1) $\quad (\Xi_N(z)^* \Xi_N(z))^{-1}_{kk} = \dfrac{\det(\Xi_{N(k)}(z)^* \Xi_{N(k)}(z))}{\det(\Xi_N(z)^* \Xi_N(z))} = \dfrac{1}{\xi_k(z)^* Q_k \xi_k(z)},$

where

$$Q_k = Q_k(z) \equiv I - \Xi_{(k)}(z)(\Xi_{(k)}(z)^* \Xi_{(k)}(z))^{-1} \Xi_{(k)}(z)^*$$

is the projection onto the orthogonal complement of the space spanned by columns $1, 2, \ldots, k-1, k+1, \ldots, N$ of $\Xi_N(z)$. The second equality in (A.1) rests on the fact that $\det(\Xi_N(z)^* \Xi_N(z)) = \det(\Xi_{N(k)}(z)^* \Xi_{N(k)}(z))$ times $(\xi_k(z)^* Q_k \xi_k(z))$, which is an instance of the Shur complement formula (see [17], page 21).

Next, since the entries of $M$, and so $X_N(z)$, are continuous random variables, each $Q_k$ is rank one with probability one. Hence, there are complex unit vectors $\gamma_k$ such that $Q_k = \gamma_k \gamma_k^*$. Summarizing, we have the bound

$$E|\operatorname{tr} \Xi_N(z)^{-1}|^p \leq N^{p-1} \sum_{k=1}^{N} E\left[\frac{1}{|\xi_k(z)^* \gamma_k|^p}\right].$$

Note that the vectors $\xi_k(z)$ and $\gamma_k$ are independent. By assumption, we also know that joint densities of the real and imaginary parts of each entry of the scaled vector $\sqrt{N}\xi_k(z)$ are uniformly bounded. From Proposition A.1, it follows that the density of the complex random variable $(\sqrt{N}\xi_k(z))^* \gamma_k$ conditional on $\gamma_k$ is bounded by a constant multiple of $N^2$. Since that constant is independent of the index $k$ and $z$, we conclude that

$$E|\operatorname{tr} \Xi_N(z)^{-1}|^p \leq N^{p-1} N^{p/2+1} \left( CN^2 \iint_{u^2+v^2 \leq 1} \frac{du\,dv}{(u^2+v^2)^{p/2}} + 1 \right)$$
$$\leq C(p) N^{3p/2+2},$$

which completes the proof. $\square$

**Acknowledgments.** This work was completed while the first author was hosted by MSRI in Berkeley, the support of which is gratefully acknowledged.

Department of Mathematics  
University of Colorado at Boulder  
UCB 395  
Boulder, Colorado 80309  
USA  
E-mail: brider@euclid.colorado.edu

Department of Mathematics  
North Carolina State University  
Box 8205  
Raleigh, North Carolina 27695  
USA  
E-mail: jack@math.ncsu.edu